\documentclass{article}

\usepackage{amssymb}
\usepackage{amsmath}
\usepackage{amsfonts}
\usepackage{amscd}
\usepackage[french]{babel}

\newtheorem{theo}{Th\'{e}or\`{e}me}
\newtheorem{defi}{D\'{e}finition}
\newtheorem{propo}{Proposition}
\newtheorem{lemm}{Lemme}

\newtheorem{coro}{Corollaire}

\usepackage[dvips]{graphicx}
\newcommand{\C}{\mathbb C}
\newcommand{\g}{\frak{g}}

\newcommand{\R}{\mathbb R}

\newcommand{\al}{alg\`ebre de Lie}

\newcommand{\h}{\frak{g} \oplus \frak{g}^*}

\begin{document}

\title{Structures complexes sur les alg\`ebres de Lie nilpotentes quasi-filiformes}
\author{Lucia Garcia-Vergnolle \footnote{Dpto. Geometr\'ia y Topolog\'ia, Facultad de Ciencias Matematicas U.C.M. Plaza de Ciencias 3, 28040 Madrid, Espagne 
lucigarcia@mat.ucm.es} et
Elisabeth~Remm \footnote{elisabeth.remm@uha.fr}\\
\\
Laboratoire de Math\'ematiques et Applications,\\
Universit\'e de Haute Alsace, Facult\'e des Sciences et,Techniques,\\
4, rue des Fr\`eres Lumi\`ere,\\
68093~Mulhouse~cedex, France.}

\date{}
\maketitle

\begin{abstract}
Le but de ce travail est  de d\'eterminer les alg\`ebres quasi-filiformes, c'est-\`a-dire dont le nilindex est 
\'egal \`a $dim(\g)-2$, qui poss\`edent
une structure complexe. Rappelons que le cas filiforme (le nilindex est \'egal \`a $dim(\g)-1$ est d\'ej\`a connu). 
\end{abstract}

{\small Mots clefs :Structures complexes. Structures complexes g\'en\'eralis\'ees}

\section{Structures complexes sur une alg\`ebre de Lie}

Soit $\g$ une alg\`ebre de Lie r\'eelle de dimension paire.
\begin{defi}
Une structure complexe sur $\g$ est un endomorphisme lin\'eaire $J$ tel que:
\begin{enumerate}
\item $J^{2}=-Id$,
\item $N(J)(X,Y)=[J(X),J(Y)]-[X,Y]-J([J(X),Y])-J([X,J(Y)])=0 \quad \forall X,Y\in \frak{g}$ \newline (condition de Nijenhuis).
\end{enumerate}
\end{defi}
Une telle structure d\'efinit une structure complexe invariante \`a gauche sur un groupe de Lie connexe
r\'eel d'alg\`ebre de Lie $\g$. Rappelons que les alg\`ebres de Lie nilpotentes munies d'une structure complexe 
sont enti\`erement d\'etermin\'ees pour les dimensions inf\'erieures ou \`egales \`a 
$6$ (voir ~\cite{Ov} et ~\cite{Sa}). Dans le cas g\'en\'eral, le seul 
r\'esultat concerne la classe des alg\`ebres de Lie filiformes. Rappelons qu'
une alg\`ebre de Lie $\frak{g}$ est dite filiforme si son nilindice est maximal, c'est-\'a-dire, 
s'il est \'egal \`a dim$(\frak{g})-1$.
\begin{propo} ~\cite{Go-Re}
Si $\frak{g}$ est une alg\`ebre de Lie filiforme de dimension paire 
alors elle n'admet pas de structures complexes.
\end{propo}
Dans ~\cite{Go-Re}, on montre dans un premier temps la non-existence de structures complexes sur 
l'alg\`ebre de Lie filiforme $\frak{L_{2n}}$ ($n\ge 2$) d\'efinie dans la base $\{X_{0},X_{1},\dots,X_{2n}\}$ par:
$$
\begin{array}{ll}
\lbrack X_{0},X_{i}\rbrack=X_{i+1}, & 1\le i \le 2n-1,\\
\lbrack X_{i},X_{j}\rbrack=0, & i,j\neq 0\\
\end{array}
$$
On g\'en\'eralise ensuite ce r\'esultat \`a toute alg\`ebre de Lie filiforme $\frak{g}$ de dimension $2n$ en remarquant que l'existence d'une structure complexe implique une d\'ecomposition
en sous-alg\`ebres $\frak{g}=\frak{g}_1 \oplus \frak{g}_2 $ o\`u $\frak{g}_1$ et $\frak{g}_2$ sont de dimension $n$. 
Une telle d\'ecomposition est impossible dans $\frak{L_{2n}}$ et donc dans toute   d\'eformation 
de  cette alg\`ebre. Comme toute  alg\`ebre de Lie filiforme de dimension $2n$ est une 
d\'eformation de $\frak{L_{2n}}$, on en d\'eduit le r\'esultat.

Cette proposition a ensuite \'etait d\'emontr\'ee dans \cite{Gu-Ca}. L'approche
 est tout \`a fait diff\'erente et repose sur la notion de structures complexes g\'en\'eralis\'ees. Nous allons
utiliser cette approche pour examiner l'existence de structures complexes sur d'autres
classes d'alg\`ebres de Lie nilpotentes.

\section{Structures complexes g\'en\'eralis\'ees sur une alg\`ebre de Lie}

\subsection{D\'efinition et lien avec les Structures Complexes}
La notion de structure complexe est d\'efinie dans le cadre g\'en\'eral des vari\'et\'es diff\'erentiables. Dans ce travail, nous allons nous int\'eresser
essentiellement au cas des structures complexes g\'en\'eralis\'ees invariantes \`a gauche sur un groupe de Lie. La d\'efinition
est alors de nature purement alg\'ebrique et s'exprime uniquement en terme d'alg\`ebres de Lie. C'est dans ce cadre que nous allons rappeler cette d\'efinition.

Soit $\g$ une alg\`ebre de Lie r\'eelle de dimension $2n$. Notons par $\g^*$ l'espace vectoriel dual qui s'identifie \`a l'espace des formes
diff\'erentielles de degr\'e $1$ invariantes \`a gauche sur un groupe de Lie connexe d'alg\`ebre de Lie $\g$. Il existe donc sur $\g^*$
une notion de diff\'erentielle ext\'erieure. Par exemple, si $\alpha \in \g^*$, alors $d\alpha \in \Lambda ^2(\g^*)$ et est donn\'ee par 
$d\alpha (X,Y)=-\alpha [X,Y]$ o\`u $[,]$ est le crochet de $\g$. Sur l'espace vectoriel $\g \oplus \g ^*$, on d\'efinit une multiplication, appel\'ee 
crochet de Courant, par:
$$
[X+\xi,Y+\eta]_c=[X,Y]+\mathcal{L}_{X}\eta-\mathcal{L}_{Y}\xi-\frac{1}{2}\textrm{d}(I_{X}{\eta}-I_{Y}\xi).
$$
pour tout $X,Y \in \g$, $\xi ,\eta \in \g^*$ et $I_{X}{\eta}$ d\'esigne le produit int\'erieur de $X$ sur $\eta$.
Cette op\'eration est antisym\'etrique et v\'erifie l'identit\'e de Jacobi (notons que dans le cadre g\'en\'eral des 
vari\'et\'es diff\'erentiables, ce crochet
se d\'efinit sur la somme du fibr\'e tangent et du fibr\'e ext\'erieur, mais ce crochet ne v\'erifie pas n\'ecessairement
l'identit\'e de Jacobi). Ainsi $\g \oplus \g ^*$, muni du crochet de Courant, est une alg\`ebre de Lie r\'eelle de dimension
$4n$. Cette alg\`ebre de Lie est une alg\`ebre de Lie quadratique. En effet 
 il existe aussi un produit scalaire donn\'e par:
$$
\langle X+\xi,Y+\eta \rangle=\frac{1}{2}(\xi(Y)+\eta(X)).
$$
\begin{defi}
Soit $\g$ une alg\`ebre de Lie r\'eelle de dimension paire $2n$. Une structure complexe g\'en\'eralis\'ee 
sur $\g$ est un endomorphisme lin\'eaire  $\mathcal{J}$ de $\g \oplus \g^*$ tel que:
\begin{enumerate}
\item $\mathcal{J}^{2}=-Id$,
\item $\mathcal{J}$ est orthogonal pour le produit scalaire $\langle \;,\; \rangle$, c'est-\'a-dire:
$$
\langle \mathcal{J}(X+\xi),\mathcal{J}(Y+\eta) \rangle=\langle X+\xi,Y+\eta \rangle \quad \forall X,Y \in \chi(M)\; \forall \xi,\eta \in \chi^{*}(M),
$$
\item Si $L$ est l'espace propre de $\mathcal{J}$ correspondant 
\`a la valeur propre $+i$, alors $L$ doit \^ etre involutif par rapport au crochet de Courant, c'est-\`a-dire $[L,L]_c\subset L$.
\end{enumerate}
\end{defi}
Remarquons que le produit scalaire est de signature $(2n,2n)$. La sous-alg\`ebre $L$ de $\h$ est un espace isotrope
$$\langle X+\xi,Y+\eta \rangle=0$$
pour tout $ X+\xi,Y+\eta \in L$. Comme il est de dimension $2n$, il est maximal isotrope. Consid\'erons sa projection sur
$\g$. On notera par $k$ la codimension de la projection de $L$ sur $\g$. Il est clair que 
$$0 \leq k\leq n.$$
\begin{defi}
Si $k$ est la codimension de la projection de $L$ sur $\g$, on dit que la structure complexe g\'en\'eralis\'ee $\mathcal{J}$ 
est de type $k$.
\end{defi}

\medskip

\noindent {\bf{Exemples}}

1. Soit $\frak{g}$ une alg\`ebre de Lie r\'eelle de dimension $2n$ munie d'une structure complexe (classique) que nous 
noterons $J$: 
$$J:\frak{g} \rightarrow \frak{g}$$
et cette application v\'erifie
$$J^2 = -Id$$
et la condition de Nijenhuis
$$N(J)(X,Y)=0$$
pour tout $X,Y \in \g.$ Cette structure permet de d\'efinir 
 une structure complexe g\'en\'eralis\'ee 
$$\mathcal{J}_{J}:\frak{g}\oplus\frak{g}^{*} \rightarrow \frak{g}\oplus\frak{g}^{*}$$
en posant
$$
\mathcal{J}_{J}(X+\xi)=-J(X)+J^{*}(\xi) \quad \forall X \in \frak{g} \; \forall \xi \in \frak{g}^{*}.
$$
Il est facile de v\'erifier que $\mathcal{J}_{J}^{2}=-Id$ et l'orthogonalit\'e de $\mathcal{J}_{J}$. 
Si on note par $T_{+}$ et $T_-$ les espaces propres de $J$ associ\'es  aux valeurs  propres $+i$ et $-i$ 
alors le $+i$-espace propre de $\mathcal{J}_{J}$ est:
$$
L=T_{-} \oplus (T_{+})^{*}
$$

On constate que l'involutivit\'e de $L$ par rapport au crochet de Courant est \'equivalente \`a ce que 
$T_{-}$ soit une sous-alg\`ebre de $\frak{g}$. Cette structure g\'en\'eralis\'ee est donc de type $n$. 

\medskip

2. Soit $\g$ une \al   r\'eelle de dimension $2n$ munie d'une forme symplectique $\omega $, c'est-\`a-dire d'une forme de degr\'e $2$
antisym\'etrique v\'erifiant
$$
\left\{
\begin{array}{l}
\omega ^n=\omega \wedge \omega \cdots \wedge \omega \neq 0 \\
d\omega (X,Y,Z)= \omega ([X,Y],Z)+\omega ([Y,Z],X)+\omega ([Z,X],Y)=0\\
\end{array}
\right.
$$
Une telle forme induit un isomorphisme, toujours not\'e $\omega$:
$$\omega : \g \longrightarrow \g ^*$$
par $\omega(X)=I_X\omega.$ 
On d\'efinit alors la structure complexe g\'en\'eralis\'ee 
$$\mathcal{J}_{\omega}:\frak{g}\oplus\frak{g}^{*} \rightarrow \frak{g}\oplus\frak{g}^{*}$$
de la fa\c{c}on suivante:
$$
\mathcal{J}_{\omega}(X+\xi)=\omega(X)+\omega^{-1}(\xi) \quad \forall X \in \frak{g} \; \forall \xi \in \frak{g}^{*}.
$$
Il s'agit d'une structure complexe g\'en\'eralis\'ee du type $0$ puisque son $+i$-espace propre est 
$$L=\{X-i\,I_{X}\omega \, : \, X \in \frak{g} \otimes  \C \}.$$

Nous avons, ci-dessus, donner les deux cas extr\^emes de structures complexes g\'en\'eralis\'ees. En g\'en\'eral, 
d'apr\`es  \cite{Gu}, toute structure complexe g\'en\'eralis\'ee du type $k$ peut s'\'ecrire comme 
une somme directe d'une structure complexe de dimension $k$ et d'une structure symplectique de dimension $2n-2k$. On en d\'eduit que toute structure de type $0$ est 
d\'efinie \` a partir d'une structure complexe sur $\g$ et que toute structure de type $n$ est donn\'ee par une forme symplectique sur $\g$.

\subsection{Approche Spinorielle}

Soit $T$ l'alg\`ebre tensorielle de $\g \oplus \g^*$ et $I$ l'ideal engendr\'e par les \'el\'ements de la forme 
$\{X+\xi \otimes X+\xi - \langle X+\xi,X+\xi \rangle . 1  \; : \; X+\xi \in \frak{g}\oplus\frak{g}^{*} \}$. L'espace quotient $C=T/I$ 
est l'alg\`ebre de Clifford de $\g \oplus \g^*$ associ\'ee au produit scalaire $\langle \, , \, \rangle$. 
Comme $C$ est une alg\`ebre associative simple, toutes les repr\'esentations irr\'eductibles de $C$ sont \'equivalentes. 
Par d\'efinition, une repr\'esentation spinorielle $\phi : C \rightarrow End_{\R}(S)$ est une repr\'esentation simple de $C$ sur $S$
et l'espace $S$ est appel\'e l'espace des spineurs.\\
\noindent Dor\'enavant, on consid\'erera $S=\wedge\frak{g}^{*}$ avec la repr\'esentation spinorielle d\'efinie par 
l'action de Clifford suivante:
$$
\begin{array}{llll}
\circ: & \frak{g}\oplus\frak{g}^{*} \times \wedge\frak{g}^{*} & \rightarrow & \wedge\frak{g}^{*}\\
 & (X+\xi,\rho) & \mapsto & (X+\xi) \circ \rho =i_{X}\rho+\xi \wedge \rho.
\end{array}
$$
Soit $\rho \in \wedge\frak{g}^{*}$  un spineur non-nul. On d\'efinit l'ensemble $L_{\rho} \subset \frak{g} \oplus \frak{g}^{*}$ par:
$$
L_{\rho}=\{X+\xi \in \frak{g}\oplus\frak{g}^{*} \; : \; (X+\xi) \circ \rho =0\}.
$$
On constate que $L_{\rho}$ est un espace isotrope. On dit que $\rho$ est un spineur pur si $L_{\rho}$ est maximal isotrope.
R\'eciproquement, si $L$ est un espace maximal isotrope, on peut consid\'erer l'ensemble $U_{L}$ des spineurs purs 
$\rho$ tels que $L=L_{\rho}$.
Dans la cas particulier o\`u $L$ est le $+i$-espace propre d'une structure complexe g\'en\'eralis\'ee, 
on peut alors prouver que l'ensemble $U_{L}$ est une droite engendr\'ee par le spineur pur:
$$
\rho=\Omega \; e^{B+i\omega}
$$
o\`u $B,\omega$ sont des $2$-formes r\'eelles et $\Omega = \theta_{1} \wedge \dots \wedge \theta_{k}$, o\` u $\theta_{1}, \dots, \theta_{k}$ 
sont des formes complexes.
De plus, on d\'eduit de  \cite{Che}  (proposition III.2.3) que $L \cap \overline{L}=\{0\}$ si et seulement si:
\begin{equation}
\label{CondBar}
\omega^{2n-2k} \wedge \Omega \wedge \overline{\Omega} \, \neq \, 0,
\end{equation}
$L$ \'etant le $+i$-espace propre de la structure complexe g\'en\'eralis\'ee.
Dans \cite{Gu}, on d\'emontre aussi que la condition d'involutivit\'e sur $L$ est \'equivalente \`a la condition d'int\'egrabilit\'e suivante:
\begin{equation}
\label{CondIn}
\exists X+\xi \in \frak{g}\oplus\frak{g}^{*} \; / \; \textrm{d} \rho = (X+\xi) \circ \rho.
\end{equation}

\subsection{Cas des alg\`ebres de Lie nilpotentes}

On consid\`ere une \al   \ r\'eelle nilpotente de dimension paire  $\frak{g}$. La suite centrale descendante est donn\'ee par:
\begin{align*}
\frak{g}^{0}  &  =\frak{g},\\
\frak{g}^{i} &=\left[  \frak{g}^{i-1},\frak{g}\right].
\end{align*}
On note  $m$ l'indice de nilpotence de $\frak{g}$. Dans $\frak{g}^{*}$, on consid\` ere la  suite croissante des sous-espaces $V_{i}$ 
o\`u les $V_{i}$ sont
des annulateurs des $\frak{g}^{i}$, c'est-\`a-dire:
$$
\left\{
\begin{array}{l}
V_{0}    =\left\{ 0 \right\} \\
V_{i}=\left\{ \varphi \in \frak{g}^{*}\backslash \; \varphi(X)=0 \, , \ \forall X \in \frak{g}^{i} \right\}. 
\end{array}
\right.
$$
Il est clair que $V_{m}=\frak{g}^{*}$. Notons la d\'efinition \'equivalente
$$V_i=\left\{ \varphi \in \frak{g}^{*}\backslash \; I_X\textrm{d}\varphi \in V_{i-1}\, , \ \forall X \in \frak{g} \right\}. $$
\begin{defi}
Soit $\alpha$ une $p$-forme sur $\g$. Le degr\'e de nilpotence de $\alpha$, not\'e par nil$(\alpha)$, est le plus petit entier 
$i$ tel que $\alpha \in \wedge^{p}V_{i}$.
\end{defi}

Supposons que $\g$ soit munie d'une structure complexe g\'en\'eralis\'ee du type $k$.
On peut ordonner les formes $\{ \theta_{1} , \dots , \theta_{k}\}$ 
par leur degr\'e de nilpotence et  les choisir de fa\c{c}on que $\{ \theta_{j}:\textrm{nil}(\theta_{j})>i\}$ 
soient lin\'eairement ind\'ependantes modulo $V_{i}$. On montre dans \cite{Gu-Ca}
qu'il existe une d\'ecomposition 
 $\Omega = \theta_{1} \wedge \dots \wedge \theta_{k}$ telle que:
\begin{description}
\item{a) }$\textrm{nil}(\theta_{i}) \leq \textrm{nil}(\theta_{j})$ si $i<j$,
\item{b) }pour chaque $i$, les formes $\{ \theta_{j}:\textrm{nil}(\theta_{j})>i\}$ sont lin\'eairement ind\'ependantes modulo $V_{i}$.
\end{description}
Une telle d\'ecomposition sera dite appropri\'ee.
\begin{theo}
Si $\g$ est une alg\` ebre de Lie nilpotente munie d'une structure complexe g\'en\'eralis\'ee, le spineur pur $\rho$ correspondant
est une forme ferm\'ee.
\end{theo}
On en d\'eduit
\begin{coro}
\label{dtheta}
Si on choisit une d\'ecomposition appropri\'ee pour $\Omega$,  alors 
\begin{description}
\item{a) }
$d\theta_{i} \in \mathcal{I}(\{ \theta_{j}: \, \textrm{nil}(\theta_{j})<\textrm{nil}(\theta_{j})\}).$
En particulier 
$$d\theta_{i} \in \mathcal{I}( \theta_{1} \dots \theta_{i-1} ).$$
\item{b) } Si $\textrm{dim} (\frac{V_{j+1}}{V_{j}})=1$ alors, ou bien il existe un $\theta_{i}$ de degr\'e de nilpotence $j$ 
ou bien il n'en n'existe pas de degr\'e $j+1$.
\end{description}
\end{coro}
{\bf Remarque.}
Suppposons qu'il existe un certain $j>0$ \`a partir duquel:
$$
\textrm{dim} (\frac{V_{i+1}}{V_{i}})=1 \quad \forall i \ge  j;
$$
S'il n'y a pas de $\theta_{i}$ de degr\'e $s\ge j$ alors, en utilisant le r\'esultat pr\'ec\'edent, 
on d\'eduit qu'il n'en n'existe pas pour tous les degr\'es  sup\'erieurs \`a $s$. 
Ceci nous permet de donner une majoration des degr\'es de nilpotence. Du corollaire ~\ref{dtheta}, on d\'eduit
 $\textrm{nil}(\theta_{1})=1$.
Si $j>1$, on a alors  $\textrm{nil}(\theta_{2}) \le j$. En effet, dans le cas contraire il n'existerait pas de 
$\theta_{i}$ de degr\'e $j$ et donc de degr\'e sup\'erieur ce qui nous m\`enerait \`a une contradiction 
car $\textrm{nil}(\theta_{2}) > j$. Par induction, on peut \'egalement d\'emontrer que $\textrm{nil}(\theta_{i}) \le j+i-2$.
Si $j=1$, par un raisonnement analogue on obtient $\textrm{nil}(\theta_{i}) \le i$.

\begin{theo}
\label{dimV}
Soit  $\frak{g}$ une alg\`{e}bre de Lie r\'eelle nilpotente de dimension $2n$ munie d'une structure complexe g\'en\'eralis\'ee du type $k>1$.
S'il existe un entier $j>0$ tel que:
$$
\textrm{dim} (\frac{V_{i+1}}{V_{i}})=1 \quad , \ \forall i \ge  j
$$
alors $k$ est major\'e par:
$$
k \, \le \, \left\{
\begin{array}{lll}
2n-\textrm{nil}(\frak{g})+j-2 & \textrm{si} & j>1\\
2n-\textrm{nil}(\frak{g}) & \textrm{si} & j=1.\\
\end{array}
\right.
$$
\end{theo}

\noindent{\it D\'emonstration.}
Supposons $j>1$. D'apr\`es la remarque pr\'ec\'edente, $\textrm{nil}(\theta_{k}) \le j+k-2$. Alors tous les 
$\theta_{1} \dots \theta_{k}$ appartiennent \`a $V_{j+k-2}$. Comme $\Omega \wedge \overline{\Omega} \neq 0$ , on a
$$
\textrm{dim} \, V_{j+k-2} \, \ge \, 2k.
$$
Par ailleurs, dim $V_{j+k-2}=2n -$dim $(\displaystyle\frac{\frak{g}^{*}}{V_{j+k-2}})$ et de plus 
$$\displaystyle\frac{\frak{g}^{*}}{V_{j+k-2}} \, \simeq \frac{V_{\textrm{nil}(\frak{g})}}{V_{\textrm{nil}(\frak{g})-1}} \oplus \dots 
\oplus \frac{V_{j+k-1}}{V_{j+k-2}},$$ donc la dimension de $V_{j+k-2}$ est \'egale \`a $2n-$nil$(\frak{g})+j+k-2$.
En rempla\c{c}ant dans l'in\'equation ci-dessus, on obtient finalement:
$$
k \, \le \, 2n-\textrm{nil}(\frak{g})+j-2.
$$
Pour $j=1$, on proc\`ede de la m\^eme fa\c{c}on en prenant $\textrm{nil}(\theta_{k}) \le k$.

\medskip

\noindent{\bf Remarque : Application au cas filiforme}. Si $\g$ est filiforme, ce th\'eor\` eme permet de retrouver le r\'esultat de
\cite{Go-Re}. En effet $m=2n-1$, $j=1$ et donc $k<2$. Il n'existe donc pas de structure de type $n$ sauf si $n=1$ mais dans ce cas l'alg\`ebre
est ab\'elienne. Nous allons maintenant regarder le cas quasi-filiforme donn\'e par $m=2n-2$.

\section{Etude des structures complexes sur les alg\`ebres quasi-filiformes}

\subsection{Classification des alg\`ebres quasi-filiformes gradu\'ees}

Soit $\frak{g}$ une alg\`{e}bre de Lie nilpotente de nilindice $m$. Elle est naturellement filtr\'{e}e par la suite centrale descendante:\\
$$
\frak{g}^{0}    =\frak{g}\supset \frak{g}^{1} \supset \frak{g}^{2} \supset\cdots \supset
\frak{g}^{k} \supset \cdots
\supset \frak{g}^{m}=\left\{ 0 \right\}.
$$
On peut alors associer une alg\`{e}bre de Lie gradu\'{e}e, not\'{e}e  $gr(\frak{g})$, et d\'{e}finie par:
$$
{\rm gr}\frak{g} \; = \; \sum_{i=1}^{m} \, \frac{\frak{g}^{i-1}}{\frak{g}^{i}} \; = \; \sum_{i=1}^{m} W_{i}
$$
dont le crochet est donn\'{e} par:
$$[X+\frak{g}^{i},Y+\frak{g}^{j}]=[X,Y]+\frak{g}^{i+j}\quad ,  \ \forall X \in \frak{g}^{i-1}\quad , \ \forall Y \in \frak{g}^{j-1}.$$
Cette alg\` ebre est appel\'ee la gradu\'ee de $\g$. Lorsque $\g$ est isomorphe \`{a} ${\rm gr}\,\frak{g}$,  
$\frak{g}$ est dite gradu\'{e}e naturellement.
On dit que $\frak{g}$ est une alg\`{e}bre de la forme $\{p_{1},\dots,p_{m}\}$ si $\dim W_{i}=p_{i}$. 
Remarquons que $gr(\frak{g})$ est de la m\^{e}me forme que $\frak{g}$.\\
Une alg\`{e}bre de Lie nilpotente est filiforme si et seulement si elle est de la forme $\{2,1,1,\dots,1\}$. On en d\'{e}duit que l'alg\`{e}bre gradu\'{e}e d'une alg\`{e}bre filiforme est aussi filiforme.

\begin{defi}
Soit $\frak{g}$ une alg\`{e}bre de Lie nilpotente , on dit que $\frak{g}$ est quasi-filiforme si son nil\'{\i}ndice $m$ est dim$\frak{g}-2$.
\end{defi}
Si  $\frak{g}$ est quasi-filiforme, il existe deux possibilit\'{e}s:
\begin{enumerate}
\item Soit $\frak{g}$ est de la forme $t_{1}=\{p_{1}=3,p_{2}=1,p_{3}=1,\dots,p_{m}=1\}$.
\item Soit $\frak{g}$ est de la forme $t_{r}=\{p_{1}=2,p_{2}=1,\dots,p_{r-1}=1,p_{r}=2,p_{r+1}=1,\dots,p_{m}=1\}$ o\`{u} $r\in \{2,\dots,m\}$.
\end{enumerate}

\begin{propo}
\label{propo}
Soit $\frak{g}$ une alg\`{e}bre de Lie quasi-filiforme gradu\'{e}e naturellement de dimension $2n$ et de la forme $t_{r}$ o\`{u} $r\in \{1,\dots,2n-2\}$. Il existe alors une base homog\`{e}ne $\{X_{0},X_{1},X_{2},\dots,X_{2n-1}\}$ de $\frak{g}$ avec  $X_{0}$ et $X_{1}$ dans $W_{1}$, $X_{i}\in W_{i}$ pour $i\in \{2,\dots,2n-2\}$ et $X_{2n-1}\in W_{r}$ dans laquelle $\frak{g}$ est une des alg\`{e}bres d\'{e}crites ci-dessous.
\begin{enumerate}
\item Si $\frak{g}$ est de la forme $t_{1}$\\
$L_{2n-1}\oplus \mathbb{R} \quad (n\ge 2)$
$$
\lbrack X_{0},X_{i}\rbrack=X_{i+1}, \quad 1 \le i \le 2n-3.
$$
\item Si $\frak{g}$ est de la forme $t_{r}$ o\`{u} $r\in \{2,\dots,2n-2\}$
\begin{enumerate}
\item $\frak{L_{2n,r}}; \quad n\ge3, \;r\: {\rm impair}, \;3\le r\le 2n-3$
$$
\begin{array}{ll}
\lbrack X_{0},X_{i} \rbrack=X_{i+1}, &  i=1,\dots,2n-3\\
\lbrack X_{i},X_{r-i} \rbrack=(-1)^{i-1}X_{2n-1}, &  i=1,\dots,\frac{r-1}{2}\\
\end{array}
$$
\item $\frak{T_{2n,2n-3}};\quad n\ge3$
$$
\begin{array}{ll}
\lbrack X_{0},X_{i} \rbrack=X_{i+1}, &  i=1,\dots,2n-4\\
\lbrack X_{0},X_{2n-1} \rbrack=X_{2n-2},\\
\lbrack X_{i},X_{2n-3-i} \rbrack=(-1)^{i-1}X_{2n-1}, &  i=1,\dots,n-2\\
\lbrack X_{i},X_{2n-2-i} \rbrack=(-1)^{i-1} (n-1-i) X_{2n-2}, &  i=1,\dots,n-2\\
\end{array}
$$
\item $\frak{N_{6,3}}$
$$
\begin{array}{ll}
\lbrack X_{0},X_{i} \rbrack=X_{i+1}, &  i=1,2,3\\
\lbrack X_{1},X_{2} \rbrack=X_{5}, &\\
\lbrack X_{1},X_{5} \rbrack=X_{4}, &\\
\end{array}
$$

\end{enumerate}
\end{enumerate}
Les crochets non \'ecrits \'etants nuls, except\'es ceux qui d\'ecoulent de l'antisym\'etrie.
\end{propo}
Pour obtenir cette classification, il suffit de reprendre celle qui a \'et\'e faite dans cas complexe \cite{Go} et \cite{Ga}.
Par exemple, si $\frak{g}$ est une alg\`ebre quasi-filiforme de la forme $t_{3}$ et de dimension $6$, 
il existe une base $\{X_{0},X_{1},\dots,X_{5}\}$ telle que
$$
\begin{array}{l}
\lbrack X_{0},X_{i}\rbrack=X_{i+1}, \, i=1,2,3,\\
\lbrack X_{1},X_{3} \rbrack=bX_{4}, \\
\lbrack X_{1},X_{2} \rbrack=bX_{3}-X_{5}, \\
\lbrack X_{5},X_{1} \rbrack=aX_{4}. \\
\end{array}
$$
Quand $a=b=0$, $\frak{g}$ est isomorphe \`a l'alg\`ebre $\frak{L_{6,3}}$. Dans le cas contraire, on consid\`ere le changement de bases
$$
Y_{0}=\alpha X_{0},\, Y_{1}=\beta X_{1}+X_{0},\, Y_{2}=\alpha \beta X_{2},\, Y_{3}=\alpha^{2} \beta X_{3},\, Y_{4}=\alpha^{3} \beta X_{4},\, Y_{5}=-\alpha \beta^{2} X_{5}
$$
avec $
\beta=\left\{
\begin{array}{lll}
-\frac{1}{b-\sqrt{|a|}} & \textrm{si} & b\neq \sqrt{|a|}\\
-\frac{1}{2\sqrt{|a|}} & \textrm{si} & b=\sqrt{|a|}\\
\end{array}
\right.
$
 et $\alpha=b\beta +1$. Les crochets sont alors donn\'es par
$$
\begin{array}{l}
\lbrack Y_{0},Y_{i}\rbrack=Y_{i+1}, \, i=1,2,3,\\
\lbrack Y_{1},Y_{3} \rbrack=Y_{4}, \\
\lbrack Y_{1},Y_{2} \rbrack=Y_{3}+Y_{5}, \\
\lbrack Y_{5},Y_{1} \rbrack=\delta Y_{4}, \, \delta=\pm 1.\\
\end{array}
$$
En faisant un deuxi\`eme changement de base, on voit que $\frak{g}$ correspond aux alg\'ebres $\frak{T_{6,3}}$ 
pour $\delta=1$ et $\frak{N_{6,3}}$ pour $\delta=-1$. Notons que dans le cas complexe, les alg\`ebres $\frak{T_{6,3}}$ et $\frak{N_{6,3}}$ 
sont isomorphes. Au del\` a de la dimension $6$, le proc\'ed\'e de construction dans le cas complexe donne la classification r\'eelle.

\begin{coro}
Soit $\frak{g}$ une alg\`{e}bre de Lie quasi-filiforme de dimension $2n$ . Il existe une base 
$\{X_{0},X_{1},X_{2},\dots,X_{2n-1}\}$ de $\frak{g}$ telle que:
\begin{enumerate}
\item Si ${\rm gr}\frak{g}\simeq \frak{L_{2n-1}}\oplus \mathbb{R} \quad (n\ge 2)$,
$$
\begin{array}{ll}
\lbrack X_{0},X_{i}\rbrack=X_{i+1}, & 1 \le i \le 2n-3,\\
\lbrack X_{i},X_{j}\rbrack=\sum_{k=i+j+1}^{2n-2}C_{i,j}^{k}X_{k}, & 1 \le i<j \le 2n-3-i,\\
\lbrack X_{i},X_{2n-1}\rbrack=\sum_{k=i+2}^{2n-2}C_{i,2n-1}^{k}X_{k}, & 1 \le i \le 2n-4,\\
\end{array}
$$
\item Si ${\rm gr}\frak{g}\simeq \frak{L_{2n,r}} \quad n\ge3, \;r\: {\rm impair}, \;3\le r\le 2n-3$,
$$
\begin{array}{ll}
\lbrack X_{0},X_{i} \rbrack=X_{i+1}, &  i=1,\dots,2n-3\\
\lbrack X_{0},X_{2n-1}\rbrack=\sum_{k=r+2}^{2n-2}C_{0,2n-1}^{k}X_{k}, &\\
\lbrack X_{i},X_{j}\rbrack=\sum_{k=i+j+1}^{2n-1}C_{i,j}^{k}X_{k}, & 1 \le i<j \le r-i-1,\\
\lbrack X_{i},X_{j}\rbrack=\sum_{k=i+j+1}^{2n-2}C_{i,j}^{k}X_{k}, & 1 \le i<j \le 2n-3-i,\; r<i+j,\\
\lbrack X_{i},X_{2n-1}\rbrack=\sum_{k=r+i+1}^{2n-2}C_{i,2n-1}^{k}X_{k}, &1 \le i \le 2n-3-r,\\
\lbrack X_{1},X_{r-1}\rbrack=X_{2n-1}, &\\
\lbrack X_{i},X_{r-i}\rbrack=(-1)^{(i-1)}X_{2n-1}+\sum_{k=r+1}^{2n-2}C_{i,r-i}^{k}X_{k}, &2 \le i \le \frac{r-1}{2},\\
\end{array}
$$
\item Si ${\rm gr}\frak{g}\simeq \frak{T_{2n,2n-3}} \quad n\ge3$,
$$
\begin{array}{ll}
\lbrack X_{0},X_{i} \rbrack=X_{i+1}, &  i=1,\dots,2n-4\\
\lbrack X_{0},X_{2n-1} \rbrack=X_{2n-2}, &\\
\lbrack X_{i},X_{j}\rbrack=\sum_{k=i+j+1}^{2n-1}C_{i,j}^{k}X_{k}, & 1 \le i<j \le 2n-4-i,\\
\lbrack X_{1},X_{2n-4} \rbrack=X_{2n-1}, &\\
\lbrack X_{i},X_{2n-3-i}\rbrack=(-1)^{(i-1)}X_{2n-1}+C_{i,2n-3-i}^{2n-2}X_{2n-2}, & 2\le i \le n-2,\\
\end{array}
$$
\item Si ${\rm gr}\frak{g}\simeq \frak{N_{6,3}}$ alors $\frak{g}\simeq \frak{N_{6,3}}$,
$$
\begin{array}{ll}
\lbrack X_{0},X_{i} \rbrack=X_{i+1}, &  i=1,2,3\\
\lbrack X_{1},X_{2} \rbrack=X_{5}, &\\
\lbrack X_{1},X_{5} \rbrack=X_{4}. &\\
\end{array}
$$
\end{enumerate}
\end{coro}
La base $\{X_{0},X_{1},X_{2},\dots,X_{2n-1}\}$ ainsi d\'{e}finie est appel\'{e}e base adapt\'{e}e de $\frak{g}$.

\subsection{Structures complexes sur les alg\`ebres de Lie quasi-filiformes}

Dans cette partie, nous allons chercher les alg\`ebres de Lie quasi-filiformes qui poss\`edent une structure complexe et donc
 une structure complexe gen\'eralis\'ee du type  $k=n$. Si $\frak{g}$ est de la forme $t_{1}$, 
le th\'eor\`eme ~\ref{dimV} implique $k=n=2$. L'alg\`ebre $\frak{g}$ est alors isomorphe \`a $\frak{L_{3}}\oplus \mathbb{R}$. 
On v\'erifie que cette alg\`ebre admet une structure complexe associ\'ee au spineur 
$$\Omega=(\omega_{0}+i\omega_{1})\wedge(\omega_{2}+i\omega_{3})$$
 $\{\omega_{0},\omega_{1},\omega_{2},\omega_{3}\}$ \'etant la base duale de la base homog\`ene $\{X_{0},X_{1},X_{2},X_{3}\}$ de la prposition ~\ref{propo}. Supposons que $\frak{g}$ soit une alg\`ebre quasi-filiforme de la forme $t_{r}$ avec $r \ge 3$. D
'apr\`es le th\'eor\`eme ~\ref{dimV}, $n=k\leq r$. 
\begin{lemm}
\label{NilQuasi}
Soit une $\frak{g}$ une alg\`ebre quasi-filiforme de la forme $t_{r}$ avec $r \ge 3$ admettant une structure complexe g\'en\'eralis\'ee 
de type $k$. On peut alors trouver des formes $\theta_{1}\dots \theta_{k}$ associ\'ees \`a la structure complexe 
generalis\'ee v\'erifiant l'une des deux conditions:
$$
\textrm{nil}(\theta_{1})=1,\,\textrm{nil}(\theta_{2})=r,\,\textrm{nil}(\theta_{3})=r+1 \dots \textrm{nil}(\theta_{k})=r+k-2
$$
ou
$$
\textrm{nil}(\theta_{1})=1,\,\textrm{nil}(\theta_{2})=r,\,\textrm{nil}(\theta_{3})=r \dots \textrm{nil}(\theta_{k})=r+k-3
$$
Dans ce dernier cas, on a $k<r$.
\end{lemm}
\noindent{\it D\'emonstration.}
Consid\'erons une d\'ecomposition appropri\'ee $\{\theta_{1}\dots \theta_{k}\}$. On d\'eduit du corollaire (~\ref{dtheta}) 
que $\textrm{nil}(\theta_{1})=1$ et $\textrm{nil}(\theta_{2})\in \{1,2,r\}$. La condition \ref{CondBar} impose alors 
$\textrm{nil}(\theta_{2})=r$ car dim$V_{1}=1$ et dim$V_{2}=3$.
Le corollaire (~\ref{dtheta}) implique alors:
$$
\textrm{nil}(\theta_{i-1}) \le \textrm{nil}(\theta_{i}) \le r+i-2 \quad i=3,\dots,k
$$
Il y a donc deux valeurs possibles pour $\textrm{nil}(\theta_{3})$:
\begin{enumerate}
\item $\textrm{nil}(\theta_{3})=r+1$\\
Supposons que $\textrm{nil}(\theta_{4})=\textrm{nil}(\theta_{3})=r+1$, les formes $\theta_{4}$ et $\theta_{3}$ 
appartiennent alors \`a $V_{r+1}$ et puisqu'elles sont ind\'ependantes modulo $V_{r}$, on a $\textrm{dim}(\frac{V_{r+1}}{V_{r}}) \ge 2$ o\`u 
$\textrm{dim}(\frac{V_{r+1}}{V_{r}}) =1$. On en d\'eduit que $\textrm{nil}(\theta_{4})=r+2$. 
Ainsi, on peut d\'emontrer que 
$$\textrm{nil}(\theta_{i})=r+i-2, \ \ \mbox{\rm pour}  \ i=3,\dots,k.$$
\item $\textrm{nil}(\theta_{3})=r$\\
De fa\c{c}on analogue, on montre que  $\textrm{nil}(\theta_{i})=r+i-3$ pour $i=3,\dots,k$.
Dans ce cas, on remarque que, si $k=r$,  le nilindice de $\theta_{r}$ est \'egal \`a $2r-3$ et alors dim$V_{2r-3} \ge 2r$. 
Ceci est impossible car dim$V_{2r-3} = 2r-1$. Ainsi $k<r$.
\end{enumerate}

\noindent{\bf Exemple.}
\label{dim6}
Consid\'erons une alg\`ebre quasi-filiforme $\frak{g}$ de dimension $6$ d\'efinie dans la base $\{X_{0},X_{1},\dots,X_{5}\}$ par:
$$
\begin{array}{ll}
\lbrack X_{0},X_{i}\rbrack =X_{i+1}, &i=1,2,3,\\
\lbrack X_{1},X_{2}\rbrack =X_{5}, &\\
\lbrack X_{1},X_{5}\rbrack =\delta X_{4}, &\delta \in \{0,1,-1\}.\\
\end{array}
$$
Supposons que $\frak{g}$ admet une structure complexe, on peut lui associer une structure complexe g\'en\'eralis\'ee du type $k=3$ 
et un spineur:
$$
\Omega=\theta_{1} \wedge \theta_{2} \wedge \theta_{3},
$$
$\theta_{1}, \theta_{2}$ et $\theta_{3}$ \'etant des formes complexes.
Notons que cette alg\`ebre est de la forme $t_{3}$ et d'apr\`es le lemme pr\'ec\'edent les nilindices correspondants sont:
$$
\textrm{nil}(\theta_{1})=1,\; \textrm{nil}(\theta_{2})=3,\; \textrm{nil}(\theta_{3})=4.
$$
Les formes complexes $\theta_{1}, \theta_{2}$ et $\theta_{3}$ peuvent donc s'\'ecrire de la fa\c{c}on suivante:
$$
\begin{array}{l}
\theta_{1}=\lambda_{0}\omega_{0}+\lambda_{1}\omega_{1},\\
\theta_{2}=\beta_{0}\omega_{0}+\beta_{1}\omega_{1}+\beta_{2}\omega_{2}+\beta_{3}\omega_{3}+\beta_{5}\omega_{5},\\
\theta_{3}=\gamma_{0}\omega_{0}+\gamma_{1}\omega_{1}+\gamma_{2}\omega_{2}+\gamma_{3}\omega_{3}+\gamma_{4}\omega_{4}+\gamma_{5}\omega_{5}\\
\end{array}
$$
o\`{u} $\lambda_{i},\beta_{i},\gamma_{i} \in \mathbb{C}$, $\gamma_{4}$ non-nul et $\beta_{3},\beta_{5}$ ne s'annulant pas simultan\'ement. 
De plus, la condition $\theta_{1} \wedge \overline{\theta_{1}} \neq 0$, est \'equivalente \`a ce que la partie imaginaire de 
$\lambda_{0} \overline{\lambda_{1}}$ soit non-nulle.
Le corollaire (~\ref{dtheta}) implique:
$$
\left\{
\begin{array}{ll}
\beta_{5}\lambda_{0}-\beta_{3}\lambda_{1} & =0\\
-\gamma_{3}\beta_{3}\lambda_{1}+\gamma_{4}\beta_{2}\lambda_{1}+\gamma_{5}\beta_{3}\lambda_{0} &=0\\
\gamma_{4}(\beta_{5}\lambda_{1}+\delta \beta_{3}\lambda_{0}) &=0\\
-\gamma_{3}\beta_{5}\lambda_{1}-\delta \gamma_{4}\beta_{2}\lambda_{0}+\gamma_{5}\beta_{5}\lambda_{0} &=0\\
\end{array}
\right.
$$
Des premi\`ere et troisi\`eme \'equations, on d\'eduit:
$$
\lambda_{1}^{2}+ \delta \lambda_{0}^{2}=0
$$
Pour $\delta=0$, ceci nous m\`ene \`a une contradiction avec $\theta_{1} \wedge \overline{\theta_{1}} \neq 0$.
Si $\delta=-1$, on obtient $\lambda_{1}=\pm \lambda_{0}$ et comme le spineur est d\'efini \`a une constante de multiplication 
pr\`es, on peut prendre $\theta_{1}=\omega_{0} \pm \omega_{1}$ ce qui contredit aussi $\theta_{1} \wedge \overline{\theta_{1}} \neq 0$.
Finalement, quand $\delta=1$, le spineur $\Omega=(\omega_{0}+i\omega_{1})\wedge(\omega_{3}+i\omega_{5})\wedge(\omega_{2}+i\omega_{4})$ 
est associ\'e \`a une structure complexe de $\frak{g}$. Ainsi l'alg\`ebre de Lie $\g$ admet une structure complexe si et seulement si
$\delta=1$.

\medskip

\begin{theo}
Soit $\frak{g}$ une alg\`ebre  de Lie r\'eelle quasi-filiforme admettant une structure complexe. Elle est alors isomorphe ou bien \`a l'alg\`ebre de dimension $4$, $\frak{L_{3}}\oplus \mathbb{R}$, ou bien \`a l'alg\`ebre de dimension $6$, $\frak{n_{6,3}}$.
\end{theo}
\noindent{\it D\'emonstration.}
Soit $\frak{g}$ une alg\`ebre  r\'eelle quasi-filiforme et de dimension $2n$ de la forme $t_{r}$ o\`u $r\in\{1,3,\dots,2n-3\}$. Supposons que $\frak{g}$ poss\`ede une structure complexe, on peut lui faire correspondre une structure complexe g\'en\'eralis\'ee du type $k=n$.\\
Pour $r=1$, nous avons vu que $\frak{g}$ est isomorphe \`a $\frak{L_{3}}\oplus \mathbb{R}$ et que cette alg\`ebre poss\`ede 
une structure complexe. Dor\'enavant, on  supposera $r\in\{3,\dots,2n-3\}$.
En appliquant le th\'eor\`eme ~\ref{dimV} et l'in\'egalit\'e:
$$
\textrm{nil}(\theta_{k})=\textrm{nil}(\theta_{n}) \leq \textrm{nil}(\frak{g})
$$
\`a chacune des possibilit\'es du lemme ~\ref{NilQuasi}, il en r\'esulte que:
\begin{enumerate}
\item Si $\textrm{nil}(\theta_{3})=r+1$ alors $\textrm{nil}(\theta_{k})=r+k-2$ et donc:
$$
n=k\le r \le n \; \Rightarrow \; r=n.
$$

\item Si $\textrm{nil}(\theta_{3})=r$ alors $\textrm{nil}(\theta_{k})=r+k-3$ et  de plus dans ce cas $k<r$, donc:
$$
n=k < r \leq n+1 \; \Rightarrow \; r=n+1.
$$

\end{enumerate}
Par ailleurs, l'alg\`ebre gradu\'ee $gr(\frak{g})$ doit \^etre isomorphe \`a l'une des alg\`ebres de la proposition ~\ref{propo}. Ainsi
 on obtient les cas suivants:
\begin{enumerate}
\item $gr(\frak{g})\sim \frak{L_{2n,r}}; \quad n\ge3, \;r\: {\rm impair}, \;3\le r\le 2n-3$.
  \begin{enumerate}
  \item Quand $\textrm{nil}(\theta_{3})=r+1$ alors $gr(\frak{g})\sim \frak{L_{2n,n}}$ avec $n\ge3\; {\rm impair}$.\\
  Notons que pour $n=3$, on retrouve l'alg\`ebre de l'exemple ~\ref{dim6} avec $\delta=0$ qui ne poss\'edait pas de structures complexes.
  Supposons $n>3$. Si $\{X_{0},X_{1},\dots,X_{2n-1}\}$ est une base adapt\'ee de $\frak{g}$,et si
 $\{\omega_{0},\omega_{1},\dots,\omega_{2n-1}\}$ est la base duale, alors
  $$
  \begin{array}{l}
  \theta_{1}=\lambda^{0}_{1}\omega_{0}+\lambda^{1}_{1}\omega_{1},\\
  \theta_{2}=\sum_{k=0}^{n}\lambda^{k}_{2}\omega_{k}+\lambda^{2n-1}_{2}\omega_{2n-1}.
  \end{array}
  $$
  Comme $\theta_{1}\wedge d\theta_{2}=0$, en regroupant les termes $\omega_{0} \wedge \omega_{1} \wedge \omega_{n-1}$, 
$\omega_{0} \wedge \omega_{2} \wedge \omega_{n-2}$ et $\omega_{0} \wedge \omega_{3} \wedge \omega_{n-3}$ dans 
$\theta_{1}\wedge d\theta_{2}$
on d\'eduit que 
$$\lambda^{n}_{2}=\lambda^{2n-1}_{2}=0.$$
Ceci est impossible car $\textrm{nil}(\theta_{2})=n$. Il n'existe pas dans ce cas de structures complexes, sauf si $n=3$.
  \item Supposons $\textrm{nil}(\theta_{3})=r$ alors $gr(\frak{g})\sim \frak{L_{2n,n+1}}$ avec $n\ge4\; {\rm pair}$. 
On peut \'ecrire $\theta_{1}$ et $\theta_{2}$ de la fa\c{c}on suivante:
  $$
  \begin{array}{l}
  \theta_{1}=\lambda^{0}_{1}\omega_{0}+\lambda^{1}_{1}\omega_{1},\\
  \theta_{2}=\sum_{k=0}^{n+1}\lambda^{k}_{2}\omega_{k}+\lambda^{2n-1}_{2}\omega_{2n-1}.
  \end{array}
  $$
  o\`u $\{\omega_{0},\omega_{1},\dots,\omega_{2n-1}\}$ est la base duale d'une base adapt\'ee de $\frak{g}$. 
Comme $\theta_{1}\wedge d\theta_{2}=0$ les coefficients correspondants aux termes  
$\omega_{0} \wedge \omega_{1} \wedge \omega_{n}$ et $\omega_{0} \wedge \omega_{2} \wedge \omega_{n-1}$, donnent:
  $$
  \left\{
  \begin{array}{l}
  \lambda_{1}^{0}\lambda_{2}^{2n-1}-\lambda_{1}^{1}\lambda_{2}^{n+1}=0\\
  \lambda_{1}^{1}\lambda_{2}^{2n-1}=0\\
  \end{array}
  \right.
  $$
  Comme $\theta_{1} \wedge \overline{\theta_{1}} \neq 0$, on d\'eduit que $\lambda_{2}^{n+1}=\lambda_{2}^{2n-1}=0$ 
ce qui contredit $\textrm{nil}(\theta_{2})=n+1$.
  \end{enumerate}
\item $gr(\frak{g})\sim \frak{T_{2n,2n-3}};\quad n\ge3$
  \begin{enumerate}
  \item Quand $\textrm{nil}(\theta_{3})=r+1$ alors $gr(\frak{g})\sim \frak{T_{6,3}}$. Dans ce cas $\frak{g}$ est isomorphe \`a l'alg\`ebre de l'exemple ~\ref{dim6} avec $\delta=-1$ qui ne poss\`ede pas de structures complexes.
  \item Quand $\textrm{nil}(\theta_{3})=r$, $gr(\frak{g})\sim \frak{T_{8,5}}$ et il existe une base adpat\'ee $\{X_{0},X_{1},\dots,X_{7}\}$de $\frak{g}$ dont les crochets v\'erifient:
  $$
  \begin{array}{ll}
  \lbrack X_{0},X_{i} \rbrack=X_{i+1}, &  i=1,\dots,4\\
  \lbrack X_{0},X_{7} \rbrack=X_{6},\\
  \lbrack X_{1},X_{i}\rbrack=\sum_{k=i+2}^{7}C_{1,i}^{k}X_{k}, & i=2,3\\
  \lbrack X_{1},X_{4} \rbrack=X_{7}, &\\
  \lbrack X_{1},X_{5} \rbrack=2X_{6}, &\\
  \lbrack X_{2},X_{4} \rbrack=-X_{6}, &\\
  \end{array}
  $$
  Dans la base duale $\{\omega_{0},\omega_{1},\dots,\omega_{7}\}$, on peut \'ecrire $\theta_{1}, \theta_{2}$ et $\theta_{3}$ 
de la fa\c{c}on suivante:
  $$
  \begin{array}{l}
  \theta_{1}=\lambda_{1}^{0}\omega_{0}+\lambda_{1}^{1}\omega_{1},\\
  \theta_{2}=\sum_{k=0}^{5}\lambda_{2}^{k}\omega_{k}+\lambda_{2}^{7}\omega_{7},\\
  \theta_{3}=\sum_{k=0}^{5}\lambda_{3}^{k}\omega_{k}+\lambda_{3}^{7}\omega_{7}.\\
  \end{array}
  $$
  D'apr\`es le corollaire ~\ref{dtheta}, $\theta_{1} \wedge d\theta_{2}=0$ et $\theta_{1} \wedge d\theta_{3}=0$.  Ainsi:
  $$
  \left\{
  \begin{array}{ll}
  \lambda_{1}^{0}\lambda_{2}^{7}-\lambda_{1}^{1}\lambda_{2}^{5} & =0\\
  \lambda_{1}^{0}\lambda_{3}^{7}-\lambda_{1}^{1}\lambda_{3}^{5} & =0\\
  \end{array}
  \right.
  $$
  En supposant $\lambda_{1}^{1}=\lambda_{2}^{7}=\lambda_{3}^{7}=1$, on obtient 
$\lambda_{1}^{0}=\lambda_{2}^{5}=\lambda_{3}^{5}$, ce qui contredit  le choix de $\theta_{2}$ et $\theta_{3}$ 
puisqu'ils sont ind\'ependants modulo $V_{4}$.
  \end{enumerate}
\item $gr(\frak{g})\sim \frak{n_{6,3}}$. $\frak{g}$ est isomorphe \`a l'alg\`ebre $\frak{n_{6,3}}$, 
l'alg\`ebre de l'exemple ~\ref{dim6} avec $\delta=1$ qui admet une structure complexe.
\end{enumerate}

\noindent Dans ~\cite{Sa}, l'alg\`ebre $\frak{n_{6,3}}$ est d\'efinie dans la base $\{X_{1},\dots,X_{6}\}$ par:
  $$
  \begin{array}{lll}
  \lbrack X_{1},X_{2} \rbrack=X_{3}, &\lbrack X_{1},X_{3} \rbrack=X_{4}, &\lbrack X_{1},X_{4} \rbrack=X_{6}, \\
  \lbrack X_{2},X_{3} \rbrack=-X_{5}, &\lbrack X_{2},X_{5} \rbrack=-X_{6}. &\\
  \end{array}
  $$
Parmi la classification de Salamon, il s'agit de la seule alg\`ebre de Lie de dimension $6$, quasi-filiforme qui poss\`ede une structure complexe.

\noindent{\bf Remerciements:}
Le premier auteur est soutenue par le projet de recherche MTM2006-09152 du Ministerio de Educaci\'{o}n y Ciencia, et remercie aussi la Fundaci\'{o}n Ram\'{o}n Areces qui finance sa bourse pr\'{e}doctorale.

\end{document}